\magnification=\magstep1

\newcount\sec \sec=0
\input Ref.macros
\input math.macros
\forwardreferencetrue
\citationgenerationtrue
\initialeqmacro
\sectionnumberstrue

\def\Isom{{\rm Isom}}
\def\Vol{{\rm Vol}}
\def\sup{{\rm sup}}

\title{Tree and Grid Factors of General Point Processes}

\author{\'Ad\'am Tim\'ar}
\abstract{We study isomorphism invariant point processes of $\R^d$ whose 
groups of symmetries are almost surely trivial. We define a 1-ended, 
locally finite tree factor on the points of the process, that is, a 
mapping of the point configuration to a graph on it
that is 
measurable and 
equivariant with the point process. This answers a question of Holroyd 
and Peres. The tree will be used to construct a factor isomorphic to 
$\Z^n$. This perhaps surprising result  (that any $d$ and $n$ works) 
solves a problem by Steve Evans. 
The construction, based on a connected clumping with $2^i$ 
vertices in each clump of the $i$'th partition, can be used to define 
various other factors.}
 
\bottomII{Primary 60G55.  Secondary
60K35. }
{Point processes, factors, random tree, random grid.}
{Research partially supported by NSF Grant 
DMS-0231224 and Hungarian National Foundation for Scientific Research 
Grants F026049 and T30074.}

\bsection{Introduction}{s.intro}

  A point process on $\R^d$ is, intuitively, a random discrete set 
of points scattered in $\R^d$. It can be thought of as a random measure 
$M$ on the Borel sets of $\R^d$ that specifies the number of points 
$M(A)$ contained in $A$ for each Borel set $A$. Given a point process $M$, 
the support of $M$ is $[M]=\{ x\in \R^d : M(\{ x\})=1 \}$, and  
points of $[M]$ are called $M$-points.
We assume throughout that the law of our point process is isometry 
invariant. 
Another property we require is that it has 
finite intensity, meaning that the expected number of $M$-points in any 
fixed bounded
Borel set $B$ is finite. Also we assume that 
an index function can be 
assigned to the set of 
$M$-points almost always, meaning that there is an injective  map from 
$[M]$  
to the real numbers 
and that it is constructed in an equivariant, measurable way. 
(A mapping $f$ from $[M]$ is equivariant with the point process if for any 
isometry $\gamma$ of 
$\R^d$, $\gamma\circ f=f\circ \gamma$.) 
 In what follows, we shall refer to this property by simply 
saying that {\it $M$ allows an index function}. 
Having an index function is equivalent to saying that the 
group of isometries of $[M]$ is trivial almost always. For this 
equivalence, see \ref b.Hol-Per/. Notice that the existence of 
an 
index function enables one to make a function (also called an index 
function) in an equivariant, measurable way that maps the pairs of $[M]$ 
to the reals 
injectively. For example, given a pair $\{x,y\}$, assume that the index of 
$x$ is $\sum_{n\in \Z} a_n 10^n$, $a_i\in 
\{0,1,...,9\}$, and that of $y$ is  $\sum_{n\in\Z} b_n 10^n$. (Thus all 
but a finite number of the $a_n$ 
and 
the $b_n$ are 0 for $n\in \Z^+$.)
Now let the index of $\{x,y\}$ be the greater of $\sum_n a_n 10^{2n} 
+\sum_k 
b_k 10^{2k+1}$ and $\sum_n b_n 10^{2n} +\sum_k 
a_k 10^{2k+1}$. 


A general example of 
point 
processes where an index function exists are 
non-equidistant processes (meaning processes that with probability 1, the 
distance between any two points is different). 
  
A {\it factor graph} or {\it factor} of $M$ is a function that maps every 
point 
configuration $[M]$ to a graph defined on it (as vertex set) and such that 
this  function 
is measurable and equivariant with the point process.
 
In 
\ref B.Fer-Lan-Tho/, it is  
shown that in dimension at most 3, a {\it translation-equivariant} 
one-ended tree 
factor of the Poisson point process
exists.  Holroyd and Peres \ref b.Hol-Per/ give a construction that 
defines a 1-ended tree 
on the 
Poisson process in an {\it isometry-equivariant} way and for any 
dimension. 
However, their proof makes use of the independence, and in the same paper 
they ask whether a one-ended tree factor can be given for any ergodic 
point 
process that almost always has only the trivial symmetry. (They give an  
example 
where this assumption about the symmetries is not satisfied, namely, the 
point process got by shifting and 
rotating $\Z^d$ in a uniform 
way. For this process, one cannot define the desired factor.) We give a 
positive answer to their question.

\procl t.treefactor  Let $M$ be a point process on $\R^d$ that allows an 
index function. Then 
there exists a locally finite one-ended tree factor on $M$.\endprocl

As shown in \ref b.Hol-Per/, a one-ended tree factor gives rise to a 
two-ended 
path factor. Briefly, define an ordering on every set of siblings 
(using the 
ordering on $[M]$ given by the index function) and then order all the 
vertices similarly to the depth-first search in 
computer 
science. 

It is also pointed out there, by a short mass-transport argument, that 
any one-ended tree that we could define 
in an 
equivariant way has to be locally finite. 
 
The same paper asks what other classes of graphs can arise as a factor of 
the Poisson process. 
For example, for what $n$ and $d$ can $\Z^n$ arise?  This question, 
due to Steve Evans, will be answered in our fourth section. There, 
we prove the following.

\procl t.gridfactor For any $d$ and $n$, a 
point process $M$ on $\R^d$ that allows an index function 
has a $\Z^n$ factor. \endprocl


\bsection{A version of the Mass-Transport Principle}{s.MTP}

The following continuum form of the Mass Transport Principle (MTP) is 
from 
\ref b.Ben-Sch/. They state it for hyperbolic spaces, but, as mentioned 
there, 
it directly generalizes to Euclidean space. 

  Call a measure $\mu$ on $\R^d\times \R^d$ {\it diagonally invariant} if 
it satisfies $$\mu (gA\times gB)= \mu (A\times B) $$
for all measurable $A, B\subset \R^d$ and $g\in \Isom (\R^d)$.

\procl l.MTP   Let $\mu$ be a nonnegative, diagonally invariant Borel 
measure
 on 
$\R^d\times \R^d$. Suppose that $\mu (A\times \R^d)< \infty$ for some 
nonempty open $A\subset \R^d$. Then 
$$\mu (B\times \R^d)=\mu (\R^d\times B)$$
for all measurable $B\subset \R^d$. Moreover, there is a constant $c$ such 
that $\mu (B\times \R^d)=c\, \Vol(B)$.\endprocl

\procl c.RadNik  Suppose that $\mu$ is a nonnegative, diagonally invariant 
Borel measure on $\R^d\times \R^d$ and that it is absolutely 
continuous with respect to 
Lebesgue measure. Let $f_{\mu}$ be its Radon-Nikod\'ym derivative: 
$\mu(A\times A')=\int_A \int_{A'} f_{\mu}(x,y)\, dx \, dy$. Then 
$\int_{\R^d} f_{\mu}(x,y)\, dx=\int_{\R^d} f_{\mu}(y,x)\ dx=c$ for almost 
every $y$, 
with the constant $c$ as in the previous lemma. \endprocl

The corollary follows from \ref l.MTP/, because if the integrals are equal 
on every Borel set then the two functions are equal almost everywhere.

We will use the lemma and its corollary in the following way. For 
convenience, we state it as a separate lemma.

\procl l.masstr  Let $T(x,y,M)$ be a nonnegative, measurable ``mass 
transport function", defined for every configuration $M$ and points 
$x,y$ of $\R^d$ 
 Suppose
$T$ is
invariant under the isometries of the space, meaning 
$T(x,y,M)=T(\gamma x,\gamma y, \gamma M)$ for any $\gamma \in 
Isom(\R^d)$. Define $f(x,y):={\bf E} T (x,y,M)$ and suppose that 
$\int_A \int_{\R^d} f(x,y)\, dx\, dy\ < \infty$ for some open $A\subset 
\R^d$.
Then $\int_{\R^d} f(x,y)\,
dx=\int_{\R^d}f(y,x)\, dx$ almost always. \endprocl

$T(x,y,M)$ is usually referred to as the amount of mass sent from $x$ to 
$y$ if the configuration is $M$. Then  $\int_{\R^d} f(x,y)\,
dx$ and $\int_{\R^d}f(y,x)\, dx$ can be thought of as the expected amount 
of mass sent into or sent out of $y$, respectively. 
 
\proof Let $\mu (A,A'):={\bf E} \int_A \int_{A'} T (x,y,M)\, dx \, dy$ and
$f(x,y):={\bf E} T (x,y,M)$. These are both isometry-invariant. 
Moreover, $f$ is the Radon-Nikod\'ym derivative of the extension of $\mu$ 
to 
${\cal B}(\R^d\times \R^d)$, since 
$$\int_A\int_{A'} f(x,y)\, dx \, dy=\int_A \int_{A'} {\bf 
E} T(x,y,M)\, dx\, dy={\bf E} 
\int_A \int_{A'} T (x,y,M)\, dx \, dy=\mu(A,A')$$ by Fubini's theorem 
for nonnegative functions. So if the assumption of the lemma about 
the existence of an open set $A$ with $\mu (A\times \R^d)< \infty$  holds, 
then the lemma and the corollary apply. 
In particular, $f_{\mu}=f$ gives $\int_{\R^d} f(x,y)\, 
dx=\int_{\R^d}f(y,x)\, dx$ almost always. \Qed

\bsection{Tree factor}{s.treefactor}

In this section we prove \ref t.treefactor/. Actually, following 
\ref b.Hol-Per/, we define a {\it 
locally finite clumping}, 
which is a sequence of coarser and coarser partitions of $[M]$, defined on 
$[M]$ in an isometry-equivariant way, and so that in every partition, all 
the classes are finite. A class in one of the partitions is 
called a {\it clump}. A clumping is {\it connected} if any 
two vertices are in the same clump in one of the  partitions (and hence 
all but in finitely many of them).  
  
As shown in \ref b.Hol-Per/, a connected locally finite clumping gives 
rise to a 
locally finite tree 
with one end. 
To construct the tree, in the first partition connect every 
vertex to the vertex of the highest index in its clump. These edges define 
a 
forest in each clump of the second partition; for each tree in this 
forest, connect the vertex of highest index in the tree to the 
 vertex of highest index
in the whole clump (but do not connect that point to itself). With these 
new edges, we 
defined a tree in each clump of the second partition, which determine a 
forest in each clump of the third partition. Continue the process this 
way. 
The graph we get after infinitely many steps is clearly a forest, 
constructed in an isometry-equivariant way. It is 
also a tree, by connectedness of the clumping. It has only one end, 
because the only path 
starting from a vertex $v$ to infinity is the one that goes through the 
vertices of greatest index
in each clump which contains $v$. 

We will need a few lemmas to construct the clumping. A subset of the 
vertex set of a graph $G$ is called {\it independent}, if no two of its 
elements are adjacent. 

\procl l.independentset Let $M$ be a point process and $G$ be a locally finite graph on 
the vertex set $[M]$, defined in an 
isometry-equivariant way. There is a subset $N$ of $[M]$ that is an 
independent set of $G$ and is defined in an equivariant, 
measurable way.\endprocl

\proof  Let $\iota$ be the index function on points of the process, and 
for a $q\in\Q$ denote 
by
  $N(q)$ the set of points $v$ such that $|\iota(v)-q| < |\iota(w)-q|$
  for all $G$-neighbors $w$ of $v$. Note that $N(q)$ is always an 
independent set (possibly empty).
  Since the union over all rational $q$ of $N(q)$ is all points of the
   process (with probability one), there exist rational numbers
   $q$ such that $\P(N(q)$ nonempty$) >0$. Call such $q$ good.
  Enumerate the rationals, and let $q([M])$ be the first good rational for 
the 
configuration $[M]$. Define $N=N(q([M]))$. \Qed

The present proof of this lemma comes from Yuval Peres, replacing the 
original, 
longer one that used a result from \ref b.Ale/.

\procl c.sparse  For all $k$, there is a nonempty subset $V_k$ of $[M]$ 
chosen in 
an
equivariant way such that the distance between any two vertices in $V_k$ 
is
at least $2^k$.\endprocl
 
\proof 
 Connect two points of $[M]$ if their distance is less than $2^k$
and
apply the lemma.  \Qed

Finally, we shall use the following simple geometric fact.

\procl l.polyhedron   Let $K\subset \R^d$ be a convex polyhedron that 
contains a ball of radius 
$r$. Then the volume of $K$ divided by the surface area of $K$ is at least 
$c\, r$, where $c>0$ is a constant depending only on $d$.\endprocl

\proof
  Connect the center $P$ of the ball to each vertex, thus subdividing
the polygon to ``pyramids", whose apices are $P$. The altitudes of the
pyramids from $P$ are at least $r$ by the hypothesis, and this gives the
claim. (The area of the bases sum up to the surface area, the volumes of 
the
pyramids to the volume of $K$.)  \Qed

%

  By \ref c.sparse/, there is a sequence $V_k$ of subsets 
of 
$[M]$, constructed in an equivariant way, such that the minimal distance 
between any two points of $V_k$ is at 
least $2^k$. Let $B_k$ be the union of the boundaries of the Voronoi cells 
on $V_k$. We show that the expected volume of the Voronoi cell 
containing some fixed point $x$ in $\R^d$ is finite.

Otherwise, define a mass-transport function $T(x,y,M)$ to be 1 if $x$ and 
$y$ are in the same Voronoi cell (say, the one corresponding to an 
$M$-point $P$) 
and if $y$ is in the ball of volume 1 around $P$. Let $T(x,y,M)$ be 0 
otherwise. Define $f(x,y)$ as ${\bf E}\, T(x,y,M)$. So $\int_{\R^d} 
f(x,y)\, dy \leq 1$. This implies also that the assumption of \ref 
l.masstr/ 
holds. However, if the expected volume of the Voronoi cell 
containing $x$ is infinite, then $\int f(y,x)\, dy=\infty$. This 
contradicts 
the lemma. Thus we proved in particular:

\procl r.boundedcells The Voronoi cells of any invariant point process are 
almost always bounded.\endprocl

  After a preparatory, intuitively clear lemma, we shall define the 
connected, locally 
finite clumping on $[M]$. Denote 
by ${\cal B}$ the set of measurable sets 
of $R^d$. 

\procl l.Oinparts  Let $O$ be a fixed point of $\R^d$.
Suppose there is an equivariant measurable partition ${\cal P}([M])={\cal 
P}$ of $\R^d$ such
that all the parts are bounded with probability 1,
and suppose that for each part $P$ in ${\cal P}$ a measurable subset of it 
is given by a measurable mapping $\phi=\phi({\cal P}, P)$. Suppose that 
$({\cal P}, \phi({\cal P},.))$ is invariant under isometries 
of 
$\R^d$.
Assume 
further, that 
for each $P\in {\cal P}([M])$, $\Vol (\phi ({\cal P},P))/ \Vol (P) 
\leq p$. Then the 
probability that $O$
lies in $\cup_{P\in {\cal P}}\phi({\cal P},P)$ is at most $p$.\endprocl

\proof  Define $T(x,y,M)$ to be $1/\Vol(P)\, (>0)$ if $y$ is in  
$P\in {\cal 
P}$ and $x$ 
is in $\phi ({\cal P},P)$.
 Let
$T(x,y,M)$ be 0 otherwise. Denote, as usual, $f(x,y):={\bf E} \,  
T(x,y,M)$. Now the expected mass sent out from $O$ is $\int_{\R^d}
f(O,y)\,
dy= {\bf P}[O\in \cup_{P\in{\cal P}}\phi({\cal P},P)]$. The 
expected mass coming 
into $O$ is
$\int_{R^d} f(y,O)\, dy\leq \sup\; 
{\Vol (\phi ({\cal P},P))}/{\Vol (P)}\leq p$.
Here
the supremum is the essential supremum over configurations of ${\cal P}$ 
of the
supremum over $P\in{\cal P}$.

So by \ref l.masstr/,  ${\bf P}[O\in \cup_{P\in {\cal P}} \phi({\cal 
P},P)]\leq p$. \Qed

Define a partition ${\cal P}_k$ of $[M]$ by saying that $x, y\in 
[M]$ 
are in the same clump of ${\cal P}_k$ iff they are in the same component 
of $\R^d\setminus \cup_{i=k}^{\infty} B_i$. This clumping 
is locally finite with probability 1 by finite intensity and the fact that 
the cells defining ${\cal P}_i$ are bounded almost always. 

\procl p.clumping The ${\cal P}_k$ define a connected, locally finite 
clumping on $[M]$.\endprocl 

\proof
  What we have to prove is that the clumping is connected. This is 
equivalent to 
saying that for any fixed ball $Q$ in $\R^d$, 
$Q$ is intersected by only a finite 
number of the $B_k$'s almost always, and so any two $M$-points inside 
$Q$ 
are in the same clump of ${\cal P}_k$ if $k$ is large enough. 
 
Denote by $\delta$ the diameter of $Q$. Now let $N_k$ be the set of 
points in 
$\R^d$ of distance less than $\delta$ 
from $B_k$, the union of the {\it thickened boundaries} of 
the 
Voronoi cells of 
$V_k$. Notice that the volume of the thickened boundary of a cell is 
bounded from 
above by $a$ times the surface area of the cell, where 
$a$ is a constant depending only on $\delta$. Here we are using that every 
cell contains a ball of radius $2^k$, by the choice of $V_k$.

 Hence by \ref l.polyhedron/, for any Voronoi 
cell on $V_k$, the volume of the cell is at least $c\, 2^k$ times as much 
as 
the volume of the thickened boundary of that cell, with some constant 
$c$ independent of 
$k$.

$Q$ is intersected by 
$B_k$ only if 
$N_k$ contains the 
center $O$ of $Q$. So it suffices to prove that for any fixed point $O$, 
the 
expected number of $N_k$'s that contain $O$ is finite. 
  
In \ref l.Oinparts/, put ${\cal P}$ to be the Voronoi cells on $V_k$ ($k$ 
fixed), and  $\phi(P)$ to be the intersection of the Voronoi cell $P$ with 
the 
thickened boundary.  
The lemma combined with \ref r.boundedcells/ says that the probability 
that $O$ is contained in 
$N_k$ is at most $2^{-k}/c$. Hence the expected number of 
$N_k$'s 
containing $O$ is at most $1/c$, and we are done. \Qed

\bsection{Grid factor}{s.gridfactor}

Once we have the 1-ended, locally finite tree factor, we can use it to 
construct a connected clumping with special properties.

\procl t.clumping There is a connected clumping $\{{\cal P}_i\}$ such that 
the 
clumps in ${\cal P}_i$ have 
size $2^i$ for each $i\geq 0$.\endprocl

\proof
 Take a 1-ended tree factor, which exists by 
\ref t.treefactor/, and denote by $T$ the actual tree given by it 
on 
the configuration $[M]$. Everything we do will be obviously equivariant 
(given 
that the tree was constructed in an equivariant way). When we have to 
decide, 
say, how to make pairs of the points of a given finite set (and which one 
to leave without a pair if there is an odd number of them), we can always 
use the index function to do this deterministically and equivariantly. We 
can say, for example, that we match the two with the highest indices, then 
the next two, etc.

Define the partition ${\cal P}_0$ to consist of singletons.

Now we define ${\cal P}_1$ in countably many steps, each step having two 
phases. 
To begin with, define two-element clumps by first forming as many pairs as 
possible in each set of leaves of $T_0 :=T$ with a common parent. Put 
these 
pairs in 
${\cal P}_1$ (so that they will be clumps of it) and delete them from 
$T_0$. We 
are left with a subtree $T_0'$ of $T_0$ such that from each set of 
sibling leaves of $T_0$, at most one is still in $T_0'$. Now, in the 
second phase, form all the pairs $\{ x, y\}$ 
such that $x$ was a leaf in $T_0$ and $y$ is its parent.

 Put these new pairs as clumps in ${\cal P}_1$, and delete 
them from $T'_0$ to get $T_1$. Observe that $T_1$ is a tree by our 
definitions. 

In the next step, do the same two phases for $T_1$ as 
 we did for 
$T_0$. Call the tree remaining at the end $T_2$, and so on. 
After 
countably many steps, all the 
vertices of $T$ are in some 
clump of ${\cal P}_1$ since the number of descendants of each vertex 
strictly 
decreases in each step (each $T_i$) until the vertex is removed from the 
actual tree. 

To define ${\cal P}_2$, identify the vertices in every pair of 
${\cal P}_1$. So 
we 
identify either connected vertices or siblings. The 
first case results in a loop; delete it. The second case results in a pair 
of 
parallel edges; delete one of the copies. The resulting graph is a 
1-ended tree 
$\hat T$, and each vertex of it represents a clump of ${\cal P}_1$.
 So the pairs on $\hat T$ defined in the same 
way as we did in the previous two paragraphs for $T$ will determine 
the ${\cal P}_2$ as we desire.

We proceed similarly to get the ${\cal P}_i$ (using ${\cal P}_{i-1}$).

Now, it is easy to see that the clumping defined is connected. That is,  
any $x,y\in V(T)$(=$[M]$) are in the same clump of ${\cal P}_k$ if $k$ is 
large enough. Indeed, we 
may 
assume that $x$ is a descendant of $y$ (otherwise choose a common 
ancestor $z$ and the bigger of the clumps containing $\{ x,z \}$ and 
$\{y,z\}$ respectively). When we defined ${\cal P}_i$ we used a tree, 
denote it by
 $\Upsilon_i$, that came 
from the tree $\Upsilon_{i-1}$ of ${\cal P}_{i-1}$ by identifying the 
pairs of a 
complete pairing of the vertices. ($\Upsilon_1$ was $T$ itself.) Thus 
every 
vertex in $\Upsilon_i$ is the result of a sequence of fusions and hence 
corresponds to $2^{i-1}$ vertices of $T$.
Denote by $v_i$  the vertex of $\Upsilon_i$ that $v$ was fused into after 
the 
sequence of $i$ identifications. Notice that if $v_i$ has at 
least one descendant, then it has 
strictly fewer descendants than $v_{i-1}$. So for $i$ large enough, $v_i$ 
has 
no descendants. For these $i$, all the descendants of $v$ in $T$ are in 
the same clump of ${\cal P}_i$. \Qed

Now we are ready to prove \ref t.gridfactor/.

\proofof t.gridfactor

Define a subgrid $K_k$ of $\Z^n$ to be the subgraph induced 
by the vertex set $V(K_k)=\{1,\ldots ,2^j \}^{n-i}\times 
\{1,\ldots ,2^{j+1}\}^i$, where $k=jn+i$, $i\in \{0,...,n-1\}$. Notice 
that 
$K_{k+1}$ arises as two 
copies of $K_k$ glued together along a ``hyperface".

Now let ${\cal P}_i$ be a clumping as in \ref t.clumping/. Use the index 
function 
to define a graph isomorphic to 
$K_i$ on each clump $C$ of ${\cal P}_i$ recursively. For the two-element 
clumps 
of ${\cal P}_1$ define $K_1$ by making the two vertices adajecent. Given 
two clumps $C_1$ and $C_2$ of ${\cal P}_i$ whose union is the clump $C$ of 
${\cal P}_{i+1}$, define a $K_{i+1}$ on $C$
by 
adding new edges to the 
union 
of the two $K_{i-1}$-graphs (defined on $C_1$ and $C_2$). If there are 
more ways to do it, the index function can be used to make it 
deterministic. 

The limiting graph $G$ is connected and clearly is a subgraph of $\Z^n$, 
because 
any finite neighbourhood of any point in it is isomorphic to a subgraph 
of some $K_i$, thus defining an embedding of $G$ in $Z^n$. 
Moreover, it is such that if $\phi$ is an embedding of $G$ to $\Z^n$, 
then  for any axis of $\Z^n$ there is at least one direction such that  
for any vertex $v$ in $\phi(G)$ the 
 infinite path ``parallel" to the axis and starting from $v$  in this 
direction is in $\phi (G)$. This 
implies that there is essentially one embedding of $G$ to $\Z^n$, meaning 
that any embedding arises from another by composing it with an isometry 
of $\Z^n$. (This is a consequence  of the fact that for any two 
subgraphs of
$Z^n$ that are both isomorphic to the graph induced by
$\{(x_1,...,x_n)\in Z^n , x_i\geq 0\}$ in $Z^n$, there is an isomorphism
between them that extends to an automorphism of $Z^n$. This claim is
intuitively obvious and one can give a proof without any difficulty.) So 
it is 
well defined 
to speak about paths in $G$ that are 
parallel to an axis - just take any embedding in $Z^n$.

Suppose now that
$G$ is 
not equal to $\Z^n$. Then for every vertex $v$ and axis $x_i$ of $\Z^n$ 
such 
that there is only one edge in $G$ incident to $v$ and parallel to the 
axis 
$x_i$, let $v$
send 
mass 1 to every point on the singly infinite path
starting 
from $v$ and parallel to $x_i$.
The expected 
mass received is 
at most $2n$, while the expected mass sent out is infinite if the limit 
graph has ``boundary points"  (points of degree less than $2n$) with 
positive probability. This contradiction with \ref l.MTP/ finishes the 
proof.
\Qed

  Let us mention that the clumping provided by \ref t.clumping/ gives rise 
to 
easy constructions of
other factors, such as 1-ended locally finite trees of arbitrary 
growth rate. For this, take the clumping with clumps of size $2^n$, where 
$n$ goes through a sequence $\{n_k\}_{k=1}^{\infty}$ of real numbers that 
go to infinity as fast as 
we wish. Then the  construction used to get a one-ended tree from a 
clumping will give us the tree that grows 
"fast".

Finally we indicate that the statements of 
 \ref t.treefactor/ and
\ref t.gridfactor/
hold also in 
the following modified
setting.  Let $G$ be any group of isometries of $R^d$, and modify the
definition of equivariance to refer to all isometries in $G$.  Then the
conclusions of the theorems hold (with the same proofs) for any
$G$-invariant process with a $G$-equivariant index function.


\medbreak
\noindent {\bf Acknowledgements.}\enspace
I thank Russell Lyons and Yuval Peres for their valuable suggestions and 
simplifications, and the anonymous referee for his carefulness.

\startbib

\bibitem[Ale]{Ale} K. S. Alexander. Percolation and minimal spanning
forests
in infinite graphs. {\it Ann. Probab.}, 23(1):87-104, 1995.
\bibitem[Ben-Sch]{Ben-Sch} I. Benjamini and O. Schramm. Percolation in
the
hyperbolic plane. {\it J. Amer. Math. Soc.}, 14(2):487-507
(electronic), 2001.
\bibitem[Fer-Lan-Tho]{Fer-Lan-Tho} P. A. Ferrari, C. Landim, and H. 
Thorisson. 
Poisson trees, succession lines and coalescing random walks. Preprint. 
\bibitem[Hol-Per]{Hol-Per} A. E. Holroyd and Y. Peres. Trees and
matchings from
point processes. {\it Elect. Comm. in Probab.}, 8:17-27 (electronic),
2003.

\endbib

\bibfile{\jobname}
\def\noop#1{\relax}
\input \jobname.bbl

\filbreak
\begingroup
\eightpoint\sc
\parindent=0pt\baselineskip=10pt

Department of Mathematics,
Indiana University,
Bloomington, IN 47405-5701
\emailwww{atimar@indiana.edu}{}
\endgroup

\bye